\documentclass[11pt]{article}
\usepackage{amssymb}
\usepackage{latexsym}
\newtheorem{thm}{Theorem}[section]
\newtheorem{lem}[thm]{Lemma}
\newtheorem{cor}[thm]{Corollary}
\newtheorem{pro}[thm]{Proposition}
\newtheorem{ex}[thm]{Example}

\newtheorem{defi}[thm]{Definition}

\newcommand{\gm }{\Gamma }
\newcommand{\lon }{\longrightarrow }
\newcommand{\be }{\begin{eqnarray*}}
\newcommand{\ee }{\end{eqnarray*}}

\input amssym.def
\input amssym

\newcommand{\pf}{\noindent{\bf Proof.}\ }
\newcommand{\qed}{\begin{flushright} $\Box$\ \ \ \ \ \
                  \end{flushright}}
\newcommand{\complex}{{\mathbb C}}
\newcommand{\reals}{{\mathbb R}}

\newcommand{\frakg}{{\mathfrak g}}
\newcommand{\frakh}{{\mathfrak h}}

\newcommand{\half}{\frac{1}{2}}
\newcommand{\third}{\frac{1}{3}}

\newcommand{\D}{{\cal D}}

\newcommand{\smalcirc}{\mbox{\tiny{$\circ $}}}

\def\description label#1{\hfil\bf[#1]\hfil}
\newcommand{\jjj}{J_{1}}
\newcommand{\jjk}{J_{2}}
\newcommand{\kk}{K_{1}}
\newcommand{\kkk}{K_{2}}
\newcommand{\deltaa}{d_{0}}

\def\sdp{\mathbin{\hbox{$\mapstochar\kern-.3333em\times$}}}
\def\pds{\mathbin{\hbox{$\times\kern-.55em\mapstochar\,$}}}

\newcommand{\wed}{\mathbin{\lower1.5pt\hbox{$\scriptstyle{\wedge}$}}}

\let\Tilde=\widetilde

\def\chigh{{\raise1.5pt\hbox{$\chi$}}}
\let\phi=\varphi

\def\til0{\Tilde{0}}

\def\dminus{\raise2pt\hbox{\vrule height1pt width 2ex}\hskip3pt}

\def\llangle{\langle\!\langle}
\def\rrangle{\rangle\!\rangle}

\def\pback#1{\mathbin{{{\lower1.2ex\hbox{$\times$}}\atop #1}}}

\def\vlra{\hbox{$\,-\!\!\!-\!\!\!-\!\!\!-\!\!\!-\!\!\!
-\!\!\!-\!\!\!-\!\!\!-\!\!\!-\!\!\!\longrightarrow\,$}}

\def\gpd{\,\lower1pt\hbox{$\longrightarrow$}\hskip-.24in\raise2pt
             \hbox{$\longrightarrow$}\,}

\def\lgpd{\,\lower1pt\hbox{$\vlra$}\hskip-1.02in\raise2pt\hbox{$\vlra$}\,}

\def\llgpd{\,\lower1pt\hbox{$\vvlra$}\hskip-1.3in\raise2pt\hbox{$\vvlra$}\,}

\begin{document}

\title{{\bf Manin triples for Lie bialgebroids}
\thanks{1991 {\em Mathematics
Subject Classification.} Primary 58F05. Secondary 17B66, 22A22, 53C99,
58H05.}}

\author{ ZHANG-JU LIU\thanks{Research partially supported by NSF of China
        and Math. Center of SEC.}\\
        Department  of Mathematics \\
        Peking University \\
        Beijing, 100871, China \\
        {\sf email: liuzj@sxx0.math.pku.edu.cn}\\
        ALAN WEINSTEIN \thanks{Research partially supported by NSF
        grant DMS93-09653.}\\
        Department of  Mathematics\\
        University of California\\
        Berkeley, CA 94720, USA\\
        {\sf email: alanw@math.berkeley.edu}\\
        and \\
        PING XU \thanks {Research partially supported by NSF
        grant DMS95-04913 and an NSF  postdoctoral fellowship.}\\
Department of Mathematics\\
The  Pennsylvania State University\\
 University Park, PA 16802, USA\\
        {\sf email: ping@math.psu.edu }}

\maketitle
\begin{abstract}
In his study of Dirac structures, a notion which includes both Poisson
structures and closed 2-forms, T. Courant introduced a bracket on the
direct sum of vector fields and 1-forms.  This bracket does not
satisfy the Jacobi identity except on certain subspaces.  In this
paper we systematize the properties of this bracket in the definition
of a {\bf Courant algebroid}.  This structure on a vector bundle
$E\rightarrow M$, consists of an antisymmetric bracket on the
sections of $E$ whose ``Jacobi anomaly'' has an explicit expression in
terms of a bundle map $E\rightarrow TM$ and a field of symmetric bilinear forms
on $E$.  When $M$ is a point, the definition reduces to that of a Lie
algebra carrying an invariant nondegenerate symmetric bilinear form.

For any Lie bialgebroid $(A,A^{*})$ over $M$ (a notion defined by
Mackenzie and Xu), there is a natural Courant algebroid structure on
$A\oplus A^{*}$ which is the Drinfel'd double of a Lie bialgebra when
$M$ is a point. Conversely, if $A$ and $A^*$ are complementary
isotropic subbundles of a Courant algebroid $E$, closed under the
bracket (such a bundle, with dimension half that of $E$, is called a
{\bf Dirac structure}), there is a natural Lie bialgebroid structure
on $(A,A^{*})$ whose double is isomorphic to $E$.  The theory of Manin
triples is thereby extended from Lie algebras to Lie algebroids.

Our work gives a new approach to bihamiltonian structures and a new
way of combining two Poisson structures to obtain a third one.  We
also take some tentative steps toward generalizing Drinfel'd's theory
of Poisson homogeneous spaces from groups to groupoids.

\end{abstract}

\section{Introduction}
The aim of this paper is to solve, in a unified way, several
mysteries which have arisen over the past few years in connection
with generalizations of the notion of Lie algebra in differential
geometry.

T. Courant \cite{Courant:1990} introduced the following antisymmetric bracket
operation on the sections of $TP\oplus T^*P$ over a manifold $P$:
\[
[X_1+ \xi_1, X_2 +\xi_2]=
[X_1,X_2]+ (L_{X_1}\xi_2-L_{X_2}\xi_1 +
d({\textstyle{\frac 12}}(\xi_1(X_2)-\xi_2(X_1))).
\]
Were it not for the last term, this would be the bracket for the
semidirect product of the Lie algebra ${\cal X}(P)$ of vector fields
with vector space $\Omega^1(P)$ of 1-forms via the Lie derivative
representation of ${\cal X}(P)$ on $\Omega^1(P)$. The last term, which
was essential for Courant's work (about which more will be said
later) causes the Jacobi identity to fail. Nevertheless, for
subbundles $E\subseteq TP\oplus T^*P$ which are maximally
isotropic for the bilinear form
$(X_1 +\xi_1,X_2+ \xi_2)_+=
\frac 12(\xi_1(X_2)+\xi_2(X_1))$, closure of $\Gamma(E)$ under the
Courant bracket implies  that the Jacobi identity {\em does}
hold on $\Gamma(E)$, because of the 
maximal isotropic
condition on $E$.  These subbundles are called {\em Dirac
structures} on $P$; the notion is a simultaneous generalization
of that of Poisson structure (when $E$ is the graph of a map
$\tilde\pi:T^*P\to TP$) and that of closed 2-form (when $E$ is the
graph of a map $\tilde\omega:TP\to T^{*}P$).

\begin{quote}
{\bf Problem 1}. Since the Jacobi identity is satisfied on certain
subspaces where $( \ , \ )_+$ vanishes, find a formula for the
{\em Jacobi anomaly}\footnote{``$+ c.p.$'' below (and henceforth)
will denote ``plus the other two terms
obtained by circular permutations of $(1,2,3).''$}
\[
[[e_1,e_2],e_3] + c.p.
\]
in terms of $( \ , \ )_+$.
\end{quote}

\bigskip
The vector space $\chi(P)\oplus \Omega^1(P)$ on which the Courant bracket
is defined is also a module over $C^\infty(P)$. Projection on the
first factor defines a map $\rho$ from $\chi(P)\oplus
\Omega^1(P)$ to derivations of $C^\infty(P)$. If one checks the
Leibniz identity which enters in the definition of a {\em Lie
algebroid} \cite{Mackenzie:book},
\[
[e_1,fe_2]=f[e_1,e_2]+(\rho(e_1)f)e_2 \,
\]
It turns out that this is not satisfied in general, but that it is
satisfied for Dirac structures. This suggests:

\begin{quote}
{\bf Problem 2}. Express the {\em Leibniz anomaly}
$[e_1,fe_2]-f[e_1,e_2]-(\rho(e_1)f)e_2$ in terms of $( \ , \ )_+$.
\end{quote}

When one is given an inner product on a Lie algebra, it is natural
to ask whether it is invariant under the adjoint representation.
Here again, a calculation turns up an {\em invariance anomaly}.

We solve problems 1 and 2 in the paper, finding an expression for
the invariance anomaly as well. The formulas obtained are so
attractive as to suggest:

\begin{quote}
{\bf Problem 3}. Generalize the Courant bracket by writing down a set
of axioms for a skew-symmetric 
bracket $\cal E\times\cal E\to\cal E$, a linear  map
$\cal E\to$  $\mbox{Der}(C^\infty(M))$, and a symmetric inner product
$\cal E\times\cal E\to C^\infty(M)$ on the space $\cal E$ of
sections of a vector bundle over $M$, and find other interesting
examples of the structure thus defined.
\end{quote}

\bigskip
Our solution of Problem 3 begins with the definition of a
structure which we call a {\em Courant algebroid}{\footnote{
We apologize to our French colleagues for possible confusion
with the nearly homonymous and somewhat less synonymous term,
``alg\`ebre de courants".}}.
Among the examples of Courant algebroids which we find are the
direct sum of any Lie bialgebroid \cite{MackenzieX:1994} and its dual, with
the bracket given by a symmetrized version of Courant's
original definition. This structure thus gives an answer as well
to:

\begin{quote}
{\bf Problem 4}. What kind of object is the double of a Lie
bialgebroid?
\end{quote}

\noindent Furthermore, within each Courant algebroid, one can consider
the maximal isotropic subbundles closed under bracket.  These more
general Dirac structures are new Lie algebroids (and sometimes Lie
bialgebroids).  Constructions in this framework applied to the Lie
bialgebroid of a Poisson manifold \cite{MackenzieX:1994} lead to new
ways of building Poisson structures and shed new light on the theory
of Poisson-Nijenhuis structures used to explicate the hamiltonian
theory of completely integrable systems \cite{K-SM:1990}. In
particular, we find a composition law for certain pairs of (possibly
degenerate) Poisson structures which generalizes the addition of
symplectic structures: namely, if $U:T^*P\to TP$ and $V:T^*P\to TP$
define Poisson structures such that $U+V$ is invertible, then
$U(U+V)^{-1}V$ again defines a Poisson structure.

When the base manifold $P$ is a point, a Lie algebroid is just a
Lie algebra. A Courant algebroid over a point turns out to be
nothing but a Lie algebra equipped with a nondegenerate ad-invariant
symmetric 2-form (sometimes called an orthogonal structure
\cite{me-re:algebres}).
(The formulas for the anomalies all involve derivatives, so they
vanish when $P$ is a point.) Such algebras and their maximal
isotropic subalgebras are the ingredients of the theory of Lie
bialgebras and Manin triples \cite{dr:quantum}. In fact, just as a
complementary pair of isotropic subalgebras in a Lie algebra with
orthogonal structure determines a Lie bialgebra, so a
complementary pair of Dirac structures in a Courant algebroid
determines a Lie bialgebroid. It is this fact, which exhibits our
theory as a generalization of the theory of Manin triples, which is
responsible for the application to Poisson-Nijenhuis pairs
mentioned above.

 We mentioned earlier that the notion of Dirac
structures was invented in order to treat in the same framework
Poisson structures, which satisfy the equation $[\pi,\pi]=0$, and
closed 2-forms, which satisfy $d\omega=0$.  One could look for a more direct
connection between these equations.

\begin{quote}
{\bf Problem 5}. What is the relation between the equations
$[\pi,\pi]=0$ and $d\omega=0$?
\end{quote}

Our solution to Problem 5 is very simple. In a Courant algebroid of
the form  $A\oplus A^*$, the double of a Lie bialgebroid, the
equation which  a skew-symmetric operator $\tilde I:A\to A^*$
must satisfy in order for its graph to be a Dirac structure
turns out to be the Maurer-Cartan equation $dI+\frac 12[I,I]=0$
for the corresponding bilinear form $I\in\Gamma(\wedge^2 A^*)$.
The structure of the original Courant algebroid
$TM\oplus TM^*$ (also viewed dually as $T^*M\oplus TM$)
is sufficiently degenerate that one of the terms in the
Maurer-Cartan equation drops out in each of the two cases.

The next problem arises from Drinfeld's study \cite{dr:poisson} of
Poisson homogeneous spaces for Poisson Lie groups.  He shows in that
paper that the
Poisson manifolds on which a Poisson Lie group $G$ acts transitively
are essentially (that is, if one deals with local rather than global
objects, as did Lie in the old days) in 1-1 correspondence with Dirac
subspaces of the double of the associated Lie bialgebra $(\frakg ,
\frakg^{*})$.  It is natural, then, to look for some kind of
homogeneous space associated to a Dirac subbundle in the double of a
Lie bialgebroid.

The object of which a Lie bialgebroid $E\lon P$ is the infinitesimal
limit is a Poisson groupoid, i.e. a Poisson manifold $\Gamma$ carrying
the structure of a groupoid with base $P$, for which the graph of
multiplication $\{(k,g,h)|k=gh\}$ is a coisotropic submanifold of
$\Gamma \times \overline{\Gamma}\times \overline{\Gamma}$.
($\overline{\Gamma }$ is $\Gamma $ with the opposite Poisson
structure.  See \cite{MackenzieX:1994}  \cite{MackenzieX:1996} 
and \cite{we:coisotropic}.)
Unlike in the case of groups, a Poisson groupoid corresponding to a
given Lie bialgebroid may exist only locally.

\begin{quote}

{\bf Problem 6}.  Define a notion of Poisson homogeneous space for a
Poisson groupoid.  Show that Dirac structures in the double of a Lie
bialgebroid $(A,A^{*})$ correspond to (local) Poisson homogeneous
spaces for the (local) Poisson groupoid $\Gamma $ associated to $(A,A^{*})$.
\end{quote}

Our solution to Problem 6 will be contained in a sequel to this 
paper \cite{LWX}.
Even if we work locally, it is somewhat complicated, since the
``homogeneous spaces'' for groupoids, which are already hard to define
in general (see \cite{br-da-ha:topological}), in this case can involve
the quotient spaces of manifolds by arbitrary foliations.

To give a flavor of our results, we mention here one example.  For the
standard Lie bialgebroid $(TM,T^{*}M)$, the associated Poisson
groupoid is the pair groupoid $M\times M$ with the zero Poisson
structure.  A Dirac structure transverse to $T^{*}M$ is the graph of a
closed 2-form $\omega $ on $M$.  The corresponding Poisson homogeneous
space for $M\times M$ is $M\times (M/{\cal F})$, where the factor
$M$ has the zero Poisson structure, and the factor $M/{\cal F}$ is the
(symplectic) Poisson manifold obtained from reduction of $M$ by the
characteristic foliation ${\cal F}$ of $\omega $.  (Of course, the
leaf space $M/{\cal F}$ might not be a manifold in any nice sense.)
Dually, our Dirac structure also defines a Poisson  homogeneous space for
the Poisson groupoid of the Lie bialgebroid $(T^{*}M,TM)$, which is
$T^{*}M$ with the operation of addition in fibres and the Poisson
structure given by the canonical 2-form.  The homogeneous space is
again $T^{*}M$, with the Poisson structure coming from the sum of the
canonical 2-form and the pullback of $\omega $ by the projection
$T^{*}M\lon M$.

\bigskip
We turn now to some problems which remain unsolved.

The only examples of Courant algebroids which we have given are the
doubles of Lie bialgebroids, i.e. those admitting a direct sum
decomposition into Dirac subbundles.  For Courant algebroids over a
point, there are many examples which are not of this type, even when
the symmetric form has signature zero, which is necessary
for such a decomposition.  For instance, we may take
the direct sum of two Lie algebras of dimension $k$ with invariant bilinear
forms, one positive definite and one negative definite.  Any isotropic
subalgebra of dimension $k$ must be the graph of an orthogonal
isomorphism from one algebra to the other.  Such an isomorphism may
not exist.  Even if it does, it might be the case that the graphs of
any two such isomorphisms must have a line in common.  (For instance,
take two copies of ${\mathfrak su}(2)$ and use the fact that every
rotation of $\reals^{3}$ has an axis.)   These examples and a further
study of Manin triples from the point of view of Lie algebras with
orthogonal structure may be found in \cite{me-re:lie}.

\begin{quote}
{\bf  Open Problem 1}.  Find interesting examples of Courant
algebroids which are not doubles of Lie bialgebroids, including
examples which admit one Dirac subbundle, but not a pair of transverse
ones.  Are there Courant algebroids which are not closely related to
finite dimensional Lie algebras, for which the bilinear form is
positive definite?
\end{quote}

In his study of quantum groups and the Knizhnik-Zamolodchikov
equation, Drinfeld \cite{dr:quasi} introduced quasi-Hopf algebras, in
which the axiom of coassociativity is weakened, and their classical
limits, the Lie quasi-bialgebras. The latter notion was studied in
depth by Kosmann-Schwarzbach \cite{ko:quasi} (see also
\cite{ba-ko:double}), who defined various structures
involving a pair of spaces in duality carrying skew symmetric brackets
whose Jacobi anomalies appear as coboundaries of other objects. Her
structures are not subsumed by ours, though, since our expression for the
Jacobi anomaly is zero when the base manifold is a point.  Jacobi
anomalies as coboundaries also appear in the theory of ``strongly
homotopy Lie algebras'' \cite{la-ma:strongly} and in recent work of
Ginzburg \cite{gi:resolution}. 
The relation of these studies to Courant
algebroids is the subject of work in progress with Dmitry Roytenberg.

\begin{quote}
{\bf Open Problem 2}. Define an interesting type of structure which
includes both the Courant algebroids and the Lie quasi-bialgebras
as special cases.
\end{quote}

At the very beginning of our study, we found that if the bracket on
a Courant algebroid is modified by the addition of a symmetric
term, many of the anomalies for the resulting asymmetric bracket
become zero.  This resembles the ``twisting'' phenomenon of Drinfeld
\cite{dr:quasi}.  

\begin{quote}
{\bf Open Problem 3}. What is the geometric meaning of such
asymmetric brackets, satisfying most of the axioms of a Lie
algebroid?
\end{quote}

The next problem is somewhat vague. The Maurer-Cartan equation
$d\alpha+\frac 12[\alpha,\alpha]$ appears as an integrability condition in the
theory of connections and plays an essential role in modern
deformation theory. (See \cite{mi:rational} and various original
sources cited therein.)

\begin{quote}
{\bf Open Problem 4}. Find geometric or deformation-theoretic
interpretations of the Maurer-Cartan equation for Dirac structures.
\end{quote}

Lie algebras, Lie algebroids and (the doubles of) Lie bialgebras
are the infinitesimal objects corresponding to Lie groups, Lie
groupoids, and (the doubles of) Poisson Lie groups respectively.
Moreover, Kosmann-Schwarzbach \cite{ko:quasi} has
studied the global objects corresponding to Lie
quasi-bialgebras, and Bangoura \cite{ba:quasi} has recently identified
the dual objects. Yet the following problem is unsolved, even for
$TM\oplus T^*M$.

\begin{quote}
{\bf Open Problem 5}. What is the global, groupoid-like object
corresponding to a Courant algebroid? In particular, what is the
double of a Poisson groupoid?
\end{quote}

A solution to Open Problem 4 might come from a solution to the next
problem.  When one passes from an object such as a Lie bialgebra or
even a Lie quasi-bialgebra to its double, the resulting object is
frequently ``nicer'' in the sense that some of the anomalies possessed
by the original object now vanish.

\begin{quote}
{\bf Open Problem 6}.  What is the double of a Courant algebroid?
\end{quote}

 Finally, we would like to
remark that many of the constructions in this paper can be carried out
at a more abstract level, either replacing the sections of a vector
bundle $E$ by a more general module over $C^\infty(P)$, as in
\cite{hu:poisson}, or in the context of local functionals on mapping
spaces as in \cite{Dorfman}  by Dorfman.

\bigskip
\noindent
{\bf  Acknowledgements.}  In addition to the funding sources mentioned
in the first footnote, we like to thank several institutions
for their hospitality while work on this project was being done: the
Isaac Newton Institute (Weinstein, Xu); the Nankai Institute for
Mathematics (Liu, Weinstein, Xu); Peking University (Xu).  Thanks go
also to Yvette Kosmann-Schwarzbach, Jiang-hua Lu,  Kirill Mackenzie
and Jim Stasheff for their helpful comments.

\section{Double of Lie bialgebroids}

\begin{defi}
\label{def:quasi-algebroid}
A {\em Courant algebroid} is a vector bundle $E\lon P$
equipped with a nondegenerate  symmetric  bilinear form
 $( \cdot , \cdot )$ on the bundle,  a  skew-symmetric
bracket $[\cdot , \cdot ]$ on $\gm (E)$
and a bundle map $\rho :E\lon TP$ such that the following
properties are  satisfied:
 \begin{enumerate}
\item For any $e_{1}, e_{2}, e_{3}\in \gm (E)$,
$[[e_{1}, e_{2}], e_{3}]+c.p.=\D  T(e_{1}, e_{2}, e_{3});$
\item  for any $e_{1}, e_{2} \in \gm (E)$,
$\rho [e_{1}, e_{2}]=[\rho e_{1}, \rho  e_{2}];$
\item  for any $e_{1}, e_{2} \in \gm (E)$ and $f\in C^{\infty} (P)$,
$[e_{1}, fe_{2}]=f[e_{1}, e_{2}]+(\rho (e_{1})f)e_{2}-
(e_{1}, e_{2})\D f ;$
\item $\rho \smalcirc \D =0$, i.e.,  for any $f, g\in C^{\infty}(P)$,
$(\D f,  \D  g)=0$;
\item for any $e, h_{1}, h_{2} \in \gm (E)$,
  $\rho (e) (h_{1}, h_{2})=([e , h_{1}]+\D (e ,h_{1}) ,
h_{2})+(h_{1}, [e , h_{2}]+\D  (e ,h_{2}) )$,
\end{enumerate}
where $T(e_{1}, e_{2}, e_{3})$ is the function on the base $P$
defined by:
\begin{equation}
\label{eq:T0}
 T(e_{1}, e_{2}, e_{3})=\third ([e_{1}, e_{2} ], e_{3})+c.p.,
\end{equation}
and
$\D :  C^{\infty}(P)\lon \gm (E)$
is  the map   defined\footnote{In this paper, 
$d_{0}$ denotes  the usual differential from
functions to 1-forms, while   $d$ will denote the
differential from functions to sections of the dual of a Lie
algebroid.}  by $\D = \half \beta^{-1}\rho^{*} d_{0}$, 
 where $\beta $ is
the isomorphism between $E$
and $E^*$ given by the bilinear form. In other words,
\begin{equation}
\label{eq:D}
(\D f , e)= \half  \rho (e) f .
\end{equation}
\end{defi}
{\bf Remark.} Introduce a twisted bracket (not antisymmetric!) on $\gm (E)$ by
$$[e ,h \tilde{]}=[e ,h]+\D (e, h).  $$
Then (iii) is equivalent to
\begin{equation}
[e_{1}, fe_{2} \tilde{]}=f [e_{1},e_{2}\tilde{]}+ (\rho (e_{1})f)e_{2};
\end{equation}
(v) is equivalent to
\begin{equation}
\rho (e) (h_{1} , h_{2})=([e ,h_{1}\tilde{]},  h_{2})+ (h_{1} ,
[e, h_{2}\tilde{]} );
\end{equation}
and (ii) and (iv) can be combined  into a single equation:
\begin{equation}
\rho [e_{1} , e_{2}\tilde{]}=[\rho e_{1}, \rho  e_{2}].
\end{equation}
It would be nice to interpret equation (i) in terms of this
twisted bracket. The geometric meaning of this
twisted bracket remains a mystery to us.

\begin{defi}
Let $E$ be a Courant algebroid. A subbundle $L$ of $E$ is
called {\em isotropic} if it  is isotropic under the
 symmetric bilinear form $( \cdot , \cdot )$. It is called  {\em integrable}
if $\gm (L)$ is closed under the bracket $[\cdot , \cdot ]$.
A {\em Dirac structure}, or {\em Dirac subbundle}, is  a subbundle $L$
which is maximally isotropic and integrable.
\end{defi}

The following proposition follows immediately from the definition.

\begin{pro}
\label{prop}
Suppose that $L$ is an integrable isotropic subbundle of
a Courant algebroid \\ $(E, \rho , [\cdot , \cdot ], ( \cdot , \cdot
)) $. 
Then $(L, \rho |_{L} , [\cdot , \cdot ])$
is a Lie algebroid.
\end{pro}

Suppose now that both $A$ and $A^{*}$ are  Lie algebroids over the base
manifold $P$,  with anchors $a$ and $a_{*}$ respectively.
Let $E$ denote their   vector bundle direct sum:
$E=A\oplus A^{*}$.
On   $E$, there exist  two natural nondegenerate
bilinear forms, one symmetric and another antisymmetric, which are
defined as follows:

\begin{equation}
\label{eq:pairing}
(X_{1}+\xi_{1} , X_{2}+\xi_{2})_{\pm}=\half (\langle \xi_{1},  X_{2}  \rangle
  \pm \langle \xi_{2} ,  X_{1}\rangle ).
\end{equation}

On $\gm (E)$, we  introduce a bracket  by

\begin{equation}
\label{eq:double}
[e_{1}, e_{2}]=([X_{1}, X_{2}]+L_{\xi_{1}}X_{2}-L_{\xi_{2}}X_{1}-d_{*}(e_{1},
 e_{2})_{-})
+ ([\xi_{1} , \xi_{2}]+L_{X_{1}}\xi_{2}-L_{X_{2}}\xi_{1} +d(e_{1}, e_{ 2})_{-}),
\end{equation}
where $e_{1}=X_{1}+\xi_{1}$ and $e_{2}=X_{2}+\xi_{2}$.\\\\

Finally, we let $\rho : E\lon TP$ be the bundle map defined by
$\rho =a +a_{*}$. That is,
\begin{equation}
\rho (X+\xi )=a(X)+a_{*} (\xi  ) , \ \ \forall X\in \gm (A) \mbox{ and }
\xi \in \gm (A^{*})
\end{equation}
It is easy to see that in this case the operator $\D$ as
defined by Equation (\ref{eq:D}) is given by
$$\D=d_{*}+d, $$
where $d_{*}: C^{\infty}(P)\lon \gm (A)$ and
$d: C^{\infty}(P)\lon \gm ( A^{*})$ are the
usual differential  operators associated to Lie algebroids \cite{MackenzieX:1994}.

 When $(A, A^{*})$ is a Lie bialgebra $(\frakg , \frakg^{*})$,
the  bracket above reduces to the  famous Lie
bracket  of Manin on the double $\frakg \oplus \frakg^{*}$.
On the other hand, if $A$ is the tangent bundle Lie algebroid
$TM$ and $A^{*}=T^{*}M $  with zero bracket,
then Equation (\ref{eq:double}) takes  the form:
$$[X_{1}+\xi_{1} , X_{2}+\xi_{2}]=[X_{1} ,X_{2}]
+\{ L_{X_{1}}\xi_{2} -L_{X_{2}}\xi_{1} +d(e_{1}, e_{ 2})_{-} \}.$$
This is the bracket first introduced by Courant \cite{Courant:1990}, then
generalized to  the  context of the formal variational calculus by Dorfman
\cite{Dorfman}.

Our work in this paper is largely motivated by an attempt to unify the two
examples above, based on the observation that Courant's bracket appears to be some
kind of ``double.''  In order to generalize Manin's construction to Lie
algebroids, it is necessary to have a compatibility condition between Lie
algebroid structures on a vector bundle and its dual.  Such a condition, providing
a definition of {\em Lie bialgebroid}, was given in
\cite{MackenzieX:1994}.  We quote here an equivalent formulation from
\cite{K-S:1994}.
\begin{defi}
A {\em Lie bialgebroid} is a dual pair $(A,A^*)$ of vector bundles equipped with
Lie algebroid structures such that the differential $d_*$ on $\Gamma(\wedge^*A)$
coming from the structure on $A^*$ is a derivation of the Schouten-type bracket
on $\Gamma(\wedge^*A)$ obtained by extension of the structure on $A$.
\end{defi}
 The
following two main theorems of this paper show that we have indeed found the
proper version of the theory of Manin triples for the Lie algebroid case.

\begin{thm}
\label{thm:main1}
If  $(A, A^{*})$ is a Lie bialgebroid, then
$E=A \oplus A^{*}$ together with $([\cdot , \cdot ], \rho , (\cdot , \cdot)_{+})$
is a Courant algebroid.
\end{thm}

Conversely, we have
\begin{thm}
\label{thm:main2}
In a Courant algebroid $(E, \rho , [\cdot , \cdot ],  ( \cdot , \cdot ))$,
suppose that  $L_{1}$ and $L_{2}$ are  Dirac subbundles  transversal
to each other, i.e., $E=L_{1}\oplus L_{2}$.
Then, $(L_{1}, L_{2} )$ is a Lie bialgebroid, where $L_{2}$  is
considered as the dual bundle of $L_{1}$  under the
pairing  $2( \cdot , \cdot )$.
\end{thm}

An immediate consequence of the theorems above
is the following duality  property of  Lie bialgebroids,
which was first proved in \cite{MackenzieX:1994}
and then by Kosmann-Schwarzbach \cite{K-S:1994} using a simpler method.

\begin{cor}
If $(A, A^{*})$  is a Lie bialgebroid,  so is $(A^{*}, A)$.
\end{cor}

\section{Jacobi anomaly}
In this section, we begin the computations leading to the proofs of our
main theorems.  Throughout this section, we assume that
$A$ is a Lie algebroid with anchor $a$ and that its dual $A^*$ is also equipped
with a Lie algebroid structure with anchor $a_{*}$.
However, we shall {\em not} assume  any compatibility conditions
between these two algebroid structures.

 For  simplicity, for  any $e_{i}=X_{i}+\xi_{i}\in \gm (E), i=1,2,3$,
 we let
$$J(e_{1}, e_{2}, e_{3})=[[e_{1}, e_{2} ], e_{3}]+c.p.$$

The main theorem of this section is the following

\begin{thm}
\label{thm:main-jacobi}
Assume that  both $(A, a)$  and $(A^{*}, a_{*})$
are Lie algebroids. Then, for $e_{i}=X_{i}+\xi_{i}\in \gm (E), i=1,2,3$,
we have
\begin{equation}
J(e_{1}, e_{2}, e_{3})=\D T( e_{1}, e_{2}, e_{3})-(\jjj +\jjk +c.p ),
\end{equation}
\end{thm}
where $$\jjj =i_{X_{3}}(d [\xi_{1} , \xi_{2} ]-L_{\xi_{1}}d \xi_{2}
+L_{\xi_{2}}d \xi_{1} )
+ i_{\xi_{3}} (d_{*} [X_{1}, X_{2}]-L_{X_{1}}d_{*}X_{2}
+L_{X_{2}}d_{*}X_{1} ), $$
and
$$\jjk = L_{d_{*}(e_{1}, e_{2})_{-}}\xi_{3}
+ [d(e_{1}, e_{2})_{-},\  \xi_{3}]
+L_{d (e_{1}, e_{2})_{-}}X_{3}+[d_{*}(e_{1}, e_{2})_{-},\ X_{3} ]. $$

We need a series of lemmas before proving this theorem.

\begin{lem}
For $e_{i}=X_{i}+\xi_{i}\in \gm (E), i=1,2,3$,
$T$ is skew-symmetric, and
\begin{equation}
\label{eq:T}
T(e_{1}, e_{2}, e_{3})=\half \{\langle [X_{1}, X_{2}], \xi_{3} \rangle
+\langle [\xi_{1}, \xi_{2}], X_{3}\rangle
+a(X_{3})(e_{1}, e_{2})_{-}-a_{*}(\xi_{3})(e_{1}, e_{2})_{-}\}+c.p.
\end{equation}

\end{lem}
\pf  The first assertion  is obvious from the definition of $T$.
For  the second one, we first have
\be
&&([e_{1}, e_{2}], e_{3})_{+}\\
&=&\half \{ \langle [X_{1}, X_{2}], \xi_{3} \rangle  + \langle L_{\xi_{1}}X_{2}, \xi_{3}\rangle
- \langle L_{\xi_{2}}X_{1} , \xi_{3}\rangle -a_{*}(\xi_{3}) (e_{1}, e_{2})_{-}\\
&&+\langle [\xi_{1}, \xi_{2}], X_{3}\rangle +\langle L_{X_{1}}\xi_{2} ,X_{3}
\rangle -\langle L_{X_{2}}\xi_{1} , X_{3}\rangle
+a(X_{3})(e_{1}, e_{2})_{-}\}\\
&=&\half \{[\langle  [X_{1}, X_{2}], \xi_{3} \rangle+\langle
[\xi_{1}, \xi_{2}], X_{3}\rangle ]+c.p. \}\\
&&+\half \{a_{*}(\xi_{1})\langle X_{2} ,\xi_{3} \rangle
-a_{*}(\xi_{2})\langle X_{1}, \xi_{3} \rangle
-a_{*}(\xi_{3})(e_{1}, e_{2})_{-}\\
&&\ \ \ \ \ \ \ \ +a(X_{1})\langle \xi_{2} ,X_{3}\rangle -a(X_{2})\langle \xi_{1} ,X_{3}\rangle +a(X_{3})(e_{1}, e_{2})_{-} \}\\
&=&\half  [\{\langle [X_{1}, X_{2}], \xi_{3} \rangle
+\langle [\xi_{1}, \xi_{2}], X_{3}\rangle
+a(X_{3})(e_{1}, e_{2})_{-}-a_{*}(\xi_{3})(e_{1}, e_{2})_{-}\}+c.p.]\\
&&+\half \rho (e_{1}) (e_{2}, e_{3})_{+}-\half  \rho (e_{2}) (e_{3}, e_{1})_{+}
\ee
Therefore, by taking the sum of its cyclic permutations, one obtains

\be
T(e_{1}, e_{2},  e_{3})&=&\third ([e_{1}, e_{2}], e_{3})_{+}+c.p.\\
&=& \half \{\langle [X_{1}, X_{2}], \xi_{3} \rangle
+\langle [\xi_{1}, \xi_{2}], X_{3}\rangle
+a(X_{3})(e_{1}, e_{2})_{-}-a_{*}(\xi_{3})(e_{1}, e_{2})_{-}\}+c.p.
\ee \qed

As a by-product, we obtain the following identity
by substituting Equation  (\ref{eq:T}) into the last step of
the computation of $([e_{1}, e_{2}], e_{3})_{+}$ in the
  proof above. This formula will be useful
later.

\begin{equation}
\label{eq:T1}
([e_{1}, e_{2}], e_{3})_{+}=T(e_{1}, e_{2},  e_{3})
+\half \rho (e_{1}) (e_{2}, e_{3})_{+}-\half  \rho (e_{2}) (e_{3}, e_{1})_{+}
\end{equation}

\begin{lem}
\label{lem:4.3}
\begin{equation}
i_{X}L_{\xi}d\eta =[\xi, L_{X}\eta ]-L_{L_{\xi}X} \eta
+[d\langle \eta  , X\rangle ,\  \xi ]
+d (a_{*}(\xi)\langle \eta , X\rangle ) -d\langle [\xi , \eta ] ,X\rangle .
\end{equation}
\end{lem}
\pf  For any $Y\in \gm (A)$,
\be
\langle i_{X}L_{\xi} d \eta , Y\rangle &=& (L_{\xi }d\eta )(X, Y)\\
&=& a_{*}(\xi )[d\eta (X, Y)]-d\eta (L_{\xi }X, Y)-
d\eta (X, L_{\xi }Y)\\
&=&a_{*}(\xi ) a(X) \langle \eta , Y\rangle -a_{*}(\xi )a(Y)\langle \eta , X\rangle
-a_{*}(\xi )\langle \eta , [X, Y]\rangle \\
&&-a(L_{\xi }X )\langle \eta , Y\rangle +a(Y)\langle \eta , L_{\xi}X\rangle
+\langle \eta , [L_{\xi} X , Y] \rangle \\
&&-a(X)  \langle \eta , L_{\xi }Y\rangle
 +a(L_{\xi }Y)\langle \eta , X\rangle +\langle \eta , [X, L_{\xi}Y]\rangle\\
&=&a_{*}(\xi )\langle L_{X}\eta , Y\rangle -a_{*}(\xi )a(Y)\langle \eta , X\rangle  +a(Y)a_{*}(\xi ) \langle \eta , X\rangle\\
&&-a(Y)\langle [\xi , \eta ], X\rangle -\langle L_{L_{\xi}X}\eta , Y\rangle
-\langle L_{X}\eta , L_{\xi }Y\rangle +\langle L_{\xi}Y, d\langle \eta ,X\rangle \rangle \\
&=&\langle [\xi , L_{X}\eta ], Y\rangle +\langle [d \langle \eta , X\rangle ,\xi ], Y\rangle \\
&&+a(Y)a_{*}(\xi )\langle \eta , X\rangle -a(Y) \langle
[\xi ,\eta ], X\rangle -\langle L_{L_{\xi }X}\eta , Y\rangle.
\ee
The lemma  follows immediately. \qed

\begin{lem}
\label{lem:T1}
\begin{equation}
\label{eq:14}
([e_{1}, e_{2}], e_{3})_{-}+c.p.=T(e_{1}, e_{2}, e_{3})+
[ \{a(X_{3})(e_{1}, e_{2})_{-}+2a_{*}(\xi_{3})(e_{1}, e_{2})_{-}
- \langle [\xi_{1}, \xi_{2}], X_{3}\rangle \} +c.p. ]
\end{equation}
\end{lem}
\pf  By definition,
$$([e_{1}, e_{2}], e_{3})_{-}+([e_{1}, e_{2}], e_{3})_{+}
=\langle [e_{1}, e_{2}]^{*}, X_{3}\rangle, $$
where  $[e_{1}, e_{2}]^{*}$ refers to  the component of $[e_{1}, e_{2}]$
in $\gm (A^{*})$.

It thus follows that
\be
&&\{([e_{1}, e_{2}], e_{3})_{-}+  c.p. \}+ 3T(e_{1}, e_{2}, e_{3})\\
&=&  \langle [e_{1}, e_{2}]^{*} , X_{3} \rangle  +c.p.\\
&=& \langle [\xi_{1},  \xi_{2}]+L_{X_{1}}\xi_{2}-L_{X_{2}}\xi_{1} +d(e_{1}, e_{ 2})_{-},
 \ \ X_{3} \rangle +c.p.\\
&=& \{\langle [\xi_{1},  \xi_{2}], X_{3}\rangle + a(X_{1})\langle \xi_{2}, X_{3}\rangle
-\langle \xi_{2} , [X_{1}, X_{3}]\rangle\\
&&+a(X_{3})(e_{1}, e_{2})_{-}-a(X_{2})\langle \xi_{1} , X_{3} \rangle
+\langle \xi_{1} , [X_{2}, X_{3}]\rangle \} +c.p.\\
&=& \{\langle [\xi_{1},  \xi_{2}], X_{3}\rangle
+2 \langle [X_{1}, X_{2}],  \xi_{3}\rangle
+3a(X_{3})(e_{1}, e_{2})_{-} \}+c.p.\\
&=&4T(e_{1}, e_{2}, e_{3})+[\{a(X_{3})(e_{1}, e_{2})_{-} +2a_{*}(\xi_{3})(e_{1},
e_{2})_{-}-\langle  [\xi_{1},  \xi_{2}], X_{3}\rangle  \}+c.p. ],
\ee
where the second from the last step follows essentially from
reorganizing cyclic permutation terms and the last
step uses Equation (\ref{eq:T}).
Equation (\ref{eq:14}) thus follows immediately. \qed
{\bf Proof of  Theorem \ref{thm:main-jacobi}}  We denote by $I_{1}$ and $I_{2}$ the
components of $J(e_{1}, e_{2},  e_{3})$ on $\gm (A^{*})$ and
$\gm (A)$ respectively.
Thus,
by definition,
\be
I_{1} &=&\{ [[\xi_{1}, \xi_{2}], \xi_{3}]+[L_{X_{1}}\xi_{2}-L_{X_{2}}\xi_{1} ,\xi_{3}]+[d(e_{1}, e_{2})_{-}, \xi_{3}]\\
&&\ \  +L_{ [X_{1}, X_{2}]+L_{\xi_{1}}X_{2}-L_{\xi_{2}}X_{1}-d_{*}(e_{1}, e_{
2})_{-}
}\xi_{3}\\
&&\ \ \ \ +L_{X_{3}}L_{X_{2}}\xi_{1}-L_{X_{3}}L_{X_{1}}\xi_{2}-L_{X_{3}}[\xi_{1}, \xi_{2}]\\
&&\ \ \ \ \ \ -d[ a(X_{3})(e_{1}, e_{ 2})_{-}] +d([e_{1}, e_{2}], e_{3})_{-} \}+ c.p.
\ee
By using the  Jacobi identity:
$[[\xi_{1}, \xi_{2}], \xi_{3}]+c.p. =0$ and the relation: $L_{[X_{1}, X_{2}]}=[L_{X_{1}}, L_{X_{2}}]$,
we can write

\begin{eqnarray}
I_{1}&=&\{[L_{X_{1}}\xi_{2}-L_{X_{2}}\xi_{1} ,\xi_{
3}]+ L_{L_{\xi_{1}}X_{2}-L_{\xi_{2}}X_{1}}\xi_{3} \nonumber \\
&&-L_{d_{*}(e_{1}, e_{2})_{-}}\xi_{3}
+ [d(e_{1}, e_{2})_{-},\  \xi_{3}]-L_{X_{3}}[\xi_{1}, \xi_{2}] \nonumber\\
&&+d([e_{1}, e_{2}], e_{3})_{-}-d(a(X_{3})(e_{1}, e_{2})_{-})\}+c.p.
\label{eq:11}
\end{eqnarray}

Now,
\be
L_{X_{3}}[\xi_{1}, \xi_{2}] &=& (di_{X_{3}}+i_{X_{3}}d )[\xi_{1}, \xi_{2}]\\
&=&d\langle X_{3}, [\xi_{1}, \xi_{2}] \rangle +
i_{X_{3}}L_{\xi_{1}}d\xi_{2}-i_{X_{3}}L_{\xi_{2}}d\xi_{1}\\
&&\ \ \ \ \ \ +i_{X_{3}}(d[\xi_{1}, \xi_{2}]-L_{\xi_{1}}d\xi_{2}+L_{\xi_{2}}d\xi_{1}).
\ee

By using Lemma \ref{lem:4.3} twice,  we can write

\begin{eqnarray}
&&L_{X_{3}}[\xi_{1}, \xi_{2}]+c.p. \nonumber \\
&=&\{ [L_{X_{1}}\xi_{2}-L_{X_{2}}\xi_{1}, \xi_{3}]+
L_{L_{\xi_{1}}X_{2}-L_{\xi_{2}}X_{1}}\xi_{3}+2[d(e_{1}, e_{2})_{-},\  \xi_{3}]
\nonumber \\
&&+2d( a_{*}(\xi_{3} )(e_{1}, e_{2})_{-}) -d\langle X_{3}, [\xi_{1}, \xi_{2}]
 \rangle
+i_{X_{3}}(d[\xi_{1}, \xi_{2}]-L_{\xi_{1}}d\xi_{2}+L_{\xi_{2}}d\xi_{1})\}+c.p.
\label{eq:12}
\end{eqnarray}

Substituting  Equation (\ref{eq:12}) into Equation (\ref{eq:11}),
we have
$$I_{1}=\{ d [ ([e_{1}, e_{2}], e_{3})_{-}-a(X_{3})(e_{1}, e_{2})_{-}
-2a_{*}(\xi_{3})(e_{1}, e_{2})_{-} +\langle [\xi_{1}, \xi_{2}], X_{3} \rangle ]
 -\kk -\kkk  \}+c.p.,
$$
where $$ \kk =i_{X_{3}}(d [\xi_{1} , \xi_{2} ]-L_{\xi_{1}}d \xi_{2}
+L_{\xi_{2}}d \xi_{1} )$$
and
$$\kkk  = L_{d_{*}(e_{1}, e_{2})_{-}}\xi_{3}
+ [d(e_{1}, e_{2})_{-},\  \xi_{3}].$$

It follows  from Lemma \ref{lem:T1} that
$$I_{1}=dT( e_{1}, e_{2}, e_{3})-\{\kk  +\kkk  +c.p.\}. $$
A similar formula for $I_{2}$ can be obtained in
a similar way.
The conclusion  follows immediately. \qed

\section{Proof of Theorem 2.5}
This section is devoted to the proof of  Theorem \ref{thm:main1}.
Throughout the section, we assume that $(A, A^{*})$ is a Lie bialgebroid  and
$E=A\oplus A^{*}$ as  in Theorem \ref{thm:main1}.
We  also let $\D: C^{\infty}(P)\lon \gm (E)$ and $\rho :E\lon TP$  be
defined as in Section 2.
To prove Theorem \ref{thm:main1}, it suffices to verify
all the five identities in Definition  \ref{def:quasi-algebroid}.
First, Equation (i) follows directly from   Theorem
\ref{thm:main-jacobi} and properties of Lie bialgebroids.
Equation (iv) is equivalent to saying that $aa_{*}^{*}$ is
skew symmetric, which is again  a property of a Lie bialgebroid
(see Proposition 3.6 in \cite{MackenzieX:1994}).
Below, we shall split the rest of the proof into several propositions.

\begin{pro}
For any $f\in C^{\infty}(P)$ and $e_{1}, e_{2}\in \gm (E)$, we
have
\begin{equation}
[e_{1}, fe_{2}]=f[e_{1}, e_{2}]+(\rho (e_{1})f)e_{2}-
(e_{1}, e_{2})_{+} \D f .
\end{equation}
\end{pro}
\pf Suppose that $e_{1}=X_{1}+\xi_{1}$ and $e_{2}=X_{2}+\xi_{2}$.
Then, we have
$$[e_{1}, fe_{2}]=[X_{1}, fX_{2}]+[X_{1}, f\xi_{2}]+[\xi_{1}, fX_{2}]
+[\xi_{1} , f\xi_{2}].$$
Here,  we have
\be
&&[X_{1}, fX_{2}] =f [X_{1},X_{2}]+(a(X_{1})f)X_{2};\\
&&[ \xi_{1} , f\xi_{2}]=f[\xi_{1}, \xi_{2}]+((a_{*}\xi_{1} )f)\xi_{2};\\
&&[X_{1}, f\xi_{2}]=f[X_{1}, \xi_{2}]+((aX_{1})f)\xi_{2}
-\half \langle X_{1}, \xi_{2} \rangle \D f ;\\
&&[\xi_{1}, fX_{2}]=f[\xi_{1},X_{2}] +((a_{*}\xi_{1} ) f)X_{2}-
\half \langle X_{2} , \xi_{1} \rangle   \D f .
\ee
The conclusion   follows  from adding up
 all the equations  above.     \qed

 \begin{pro}
\label{pro:hom}
For any $e_{1}, e_{2}\in \gm (E)$, we have
$$\rho [e_{1}, e_{2} ]=[\rho e_{1}, \rho e_{2}]. $$
\end{pro}

We need a lemma first.

\begin{lem}
 \label{lem:4.2}
If $(A, A^{*})$ is  a Lie bialgebroid  with anchors $(a, a_{*})$,
then for any $X\in \gm (A) $ and $\xi \in \gm (A^{*})$,
$$ [a(X), a_{*}(\xi )] =
a_{*}(L_{X}\xi )-a (L_{\xi }X) + a a_{*}^{*}\deltaa \langle \xi , X \rangle . $$
\end{lem}
\pf  For any $f\in C^{\infty}(P)$,
\be
&&(a a_{*}^{*}\deltaa \langle \xi , X \rangle ) f\\
&=& \langle d_{*}\langle \xi , X \rangle , df \rangle\\
&=&L_{df}\langle \xi , X \rangle\\
&=&\langle L_{df}\xi, X \rangle +\langle \xi, L_{df} X \rangle\\
&=&  -\langle  L_{\xi} df , X\rangle + \langle \xi, [X, d_{*}f]\rangle\\
&=& -a_{*}(\xi )\langle df  , X\rangle +\langle  df , L_{\xi} X \rangle
+ a(X)a_{*}(\xi )f-\langle L_{X}\xi , d_{*}f\rangle\\
&=&[a(X) , a_{*}(\xi )]f-a_{*}(L_{X}\xi )f+ a(L_{\xi} X )f,
\ee
where in the fourth equality we have used the fact that $L_{df}X=[X , d_{*}f]$,
a property of  a general Lie bialgebroid (see Proposition 3.4
 of \cite{MackenzieX:1994}).
\qed
  {\bf Proof of Proposition \ref{pro:hom}} Let
$e_{1}=X_{1}+\xi_{1}$ and $e_{2}=X_{2}+\xi_{2}$.
\be
\rho [e_{1}, e_{2}]&=&
a\{ [X_{1}, X_{2}]+L_{\xi_{1}}X_{2}-L_{\xi_{2}}X_{1}-
d_{*}(e_{1},  e_{2} )_{-}\}\\
&&+a_{*}\{[\xi_{1} ,\xi_{2} ]+L_{X_{1}}\xi_{2}-L_{X_{2}}\xi_{1} +
d(e_{1}, e_{2})_{-}\}\\
&=&a[X_{1}, X_{2}]+a(L_{\xi_{1}}X_{2})-a(L_{\xi_{2}}X_{1})-
\half a a_{*}^{*}\deltaa (\langle \xi_{1},  X_{2}\rangle- \langle\xi_{2},
 X_{1} \rangle )\\
&&+a_{*}[ \xi_{1} , \xi_{2} ]+a_{*}(L_{X_{1}}\xi_{2})-a_{*}(L_{X_{2}}\xi_{1})
+\half a_{*}a^{*} \deltaa (\langle \xi_{1},  X_{2}\rangle -\langle \xi_{2},
 X_{1} \rangle  )\\
&=&a[X_{1}, X_{2}]+[a(L_{\xi_{1}}X_{2})-a_{*}(L_{X_{2}}\xi_{1})
-aa_{*}^{*} \deltaa \langle \xi_{1},  X_{2} \rangle ]\\
&&-[ a(L_{\xi_{2}}X_{1})-a_{*}(L_{X_{1}}\xi_{2})-aa_{*}^{*}\deltaa
\langle \xi_{2} , X_{1} \rangle  ]
+a_{*}[\xi_{1}, \xi_{2}]\\
&=&[aX_{1}, aX_{2}]+[a_{*}\xi_{1}, a_{*}\xi_{2}]
+[aX_{1}, a_{*}\xi_{2}]+[a_{*}\xi_{1}, aX_{2}]\\
&=&[\rho (e_{1}) , \rho (e_{2})],
\ee
where in  the third equality  we have used
 the skew-symmetry  of  the operator $aa_{*}^{*}$, and
the second from the last follows from Lemma \ref{lem:4.2}.
\qed

\begin{pro}
For any $e, h_{1}, h_{2}\in \gm (E)$, we  have
\begin{equation}
\label{eq:17}
\rho (e) (h_{1}, h_{2})_{+}=([e , h_{1}]+\D (e ,h_{1})_{+} , \ h_{2})_{+}+
(h_{1}, \ [e , h_{2}]+\D  (e ,h_{2})_{+} )_{+}
\end{equation}
\end{pro}
\pf According to Equation (\ref{eq:T1}),
$$([e ,h_{1}], h_{2})_{+}=T(e, h_{1}, h_{2})+\half \rho(e) (h_{1}, h_{2})_{+}
-\half \rho (h_{1})(e, h_{2})_{+}$$
and
$$(h_{1}, [e, h_{2}])_{+}=T(e,   h_{2}, h_{1})+\half \rho(e) (h_{2}, h_{1})_{+}
-\half \rho ( h_{2}) (e, h_{1})_{+}.$$
 By adding these two equations, we obtain
Equation (\ref{eq:17})   immediately
since  $T(e, h_{1}, h_{2})$ is skew-symmetric with respect to
$h_{1}$ and $h_{2}$.
\qed

\section{Proof of Theorem 2.6}

This section is devoted to the proof of  Theorem \ref{thm:main2}.
We denote  sections of $L_{1}$  by letters $X, Y$, and
sections of $L_{2}$ by  $\xi , \eta $ etc.. For  any
$X\in \gm (L_{1})$ and
$\xi \in \gm (L_{2})$,  we define their pairing by

\begin{equation}
\langle \xi , X \rangle =2( \xi , X ).
\end{equation}

Since $(\cdot , \cdot )$ is nondegenerate,
 $L_{2}$ can be  considered as the dual bundle of
$L_{1}$ under this pairing. Moreover, the symmetric
bilinear form  $(\cdot , \cdot )_{+} $ on $E$   defined by Equation (\ref{eq:pairing})
coincides with the original one.

By Proposition \ref{prop}, both $L_{1}$ and $L_{2}$ are Lie algebroids, and
their anchors are given by $a=\rho |_{L_{1}}$ and $a_{*}=\rho |_{L_{2}}$
respectively.
We shall use $d: \gm (\wedge^{*}L_{2})\lon \gm (\wedge^{*+1} L_{2}     )
$ and $d_{*}: \gm (\wedge^{*}L_{1})\lon \gm (\wedge^{*+1} L_{1}     )$
to denote their induced de-Rham differentials as usual.

Equation (v) in Definition \ref{def:quasi-algebroid} implies immediately
 that the  bracket between $ X\in \gm (L_{1}) $ and $ \xi \in \gm (L_{2}) $
is
given by
\begin{equation}
[X ,\xi ]= (-L_{\xi } X +\half d_{*} \langle  \xi ,  X \rangle )
+( L_{X}\xi -\half d\langle  \xi , X \rangle  )
\end{equation}
Thus we have

\begin{pro}
Under the decomposition $E=L_{1}\oplus L_{2}$,  for sections $e_{i}\in
\gm (E)$, $i=1, 2$ if
we write $e_{i}=X_{i}+\xi_{i}$,
then the bracket $[e_{1}, e_{2}]$ is given by
Equation (\ref{eq:double}).
\end{pro}

Before proving Theorem \ref{thm:main2}, we need the following lemma.
\begin{lem}
\label{lem:exchange}
Under the assumption of Theorem \ref{thm:main2}
we have
\be
&&L_{d_{*}f}\xi =-[df , \xi]; \\
&&L_{df}X=-[d_{*}f , X],
\ee
for any $f\in C^{\infty}(P)$, $X\in \gm (L_{1})$ and
$\xi\in \gm (L_{2})$.
\end{lem}
\pf Clearly,   Equation (iv) in Definition
\ref{def:quasi-algebroid}
 yields that $a\smalcirc d_{*}=-a_{*} \smalcirc d$.
Therefore,
\begin{eqnarray}
[a_{*} \xi , a X]&=&[\rho \xi , \rho X ] \nonumber\\
&=&\rho [\xi , X] \nonumber \\
&=&\rho ( L_{\xi } X -\half d_{*} \langle \xi  , X \rangle -L_{X}\xi +\half
 d\langle \xi , X \rangle ) \nonumber\\
&= &a( L_{\xi } X -\half d_{*} \langle  \xi  , X\rangle )+ a_{*}(-L_{X}\xi +
\half d\langle \xi , X \rangle) \nonumber\\
&=&a( L_{\xi } X)-a_{*}(L_{X}\xi )+
 ( a_{*} d ) \langle  \xi , X \rangle .  \label{eq:15}
\end{eqnarray}

On the other hand,
\begin{eqnarray}
((a_{*}d )\langle  \xi , X \rangle ) f &=&
(a (d_{*}f))\langle  \xi ,X  \rangle  \nonumber \\
&=&\langle L_{d_{*}f}\xi , X\rangle +\langle \xi , [d_{*}f , X]\rangle  \nonumber \\
&=&\langle [\xi , df], X \rangle -\langle \xi , L_{X}d_{*}f\rangle
+\langle  L_{d_{*}f}\xi +[df , \xi ], X\rangle \nonumber \\
&=&a_{*}(\xi )a (X)f- \langle df , L_{\xi }X\rangle -
a(X)a_{*}(\xi )f +\langle L_{X}\xi , d_{*}f\rangle
+\langle  L_{d_{*}f}\xi +[df , \xi ],  X\rangle \nonumber \\
&=&[a_{*}(\xi ), a(X)]f-a(L_{\xi }X)f+a_{*} (L_{X}\xi )f
+\langle  L_{d_{*}f}\xi +[df , \xi ], X\rangle . \label{eq:16}
\end{eqnarray}

Comparing  Equation (\ref{eq:15}) with ( \ref{eq:16}),
we obtain
$$\langle L_{d_{*}f}\xi +[df , \xi ] , X\rangle
=0.$$
Therefore, $L_{d_{*}f}\xi  =-[df , \xi ]$. The other equation
can be proved similarly.
\qed
{\bf Proof of Theorem \ref{thm:main2} }
It follows from Theorem \ref{thm:main-jacobi} that
$\jjj +\jjk +c.p. =0$, for any $e_{1}, e_{2}$ and $e_{3}\in \gm (E)$.
Using  Lemma \ref{lem:exchange},  we have $\jjj + c.p. =0$.
In particular, if we take $e_{1}=X_{1}, \ e_{2}=X_{2}$ and
$e_{3}=\xi_{3}$, we  obtain that $i_{\xi_{3}} (d_{*} [X_{1}, X_{2}]-L_{X_{1}}d_{*}X_{2}
+L_{X_{2}}d_{*}X_{1} )=0$,
which implies the compatibility condition between
$A$ and $A^{*}$. \qed

\section{Hamiltonian operators}
Throughout this section, we will assume that $(A, A^{*})$ is a Lie bialgebroid.
Suppose that $H: A^{*}\lon A$ is a bundle map.  We denote by
$A_{H}$ the graph of $H$, considered as a subbundle of $E=A\oplus A^{*}$.
I.e., $A_{H}=\{ H\xi +\xi |\xi\in A^{*}\}$.

\begin{thm}
$A_{H}$ is a Dirac subbundle iff $H$ is skew-symmetric and
satisfies the following Maurer-Cartan type equation:
\begin{equation}
\label{eq:Maure}
d_{*}H+\half [H, H]=0.
\end{equation}
In this equation, $H$ is considered as a section of $\wedge^{2}A$.
\end{thm}
In the sequel, we shall use the same symbol to denote  a section
of $\wedge^{2}A$ and its induced bundle map if
no confusion is caused. \\
\pf It is easy to see that $A_{H}$ is isotropic
iff  $H$ is skew-symmetric.
For any $\xi , \eta \in \gm (A^{*})$,
let
\begin{equation}
[\xi , \eta ]_{H}=L_{H\xi}\eta -L_{H\eta }\xi +d \langle \xi , H\eta \rangle.
\end{equation}
Since
\be
&&[H\xi , \eta ]=L_{H\xi }\eta -L_{\eta }H\xi - \half (d-d_{*})\langle \eta ,
 H\xi  \rangle \ \ \mbox{and} \\
&&[\xi , H\eta ]=L_{\xi}H\eta -L_{H\eta }\xi +\half (d-d_{*})\langle  \xi ,
H\eta \rangle ,
\ee
then,
$$[H\xi , \eta ]+ [\xi , H\eta ]=[\xi , \eta ]_{H}+L_{\xi}H\eta -L_{\eta}H\xi
+d_{*}\langle \eta , H\xi  \rangle  .$$

On the other hand, we have the following
 formula (see \cite{K-SM:1990}):

\begin{equation}
[H\xi , H\eta ]=H[\xi , \eta ]_{H}-\half [H , H](\xi ,\eta ).
\end{equation}
Therefore, we have,
\be
[H\xi +\xi , H\eta +\eta ]&=&[H\xi , H\eta ]+[\xi , H\eta ]+[H\xi , \eta ]+
[\xi , \eta ]\\
&=&(L_{\xi}H\eta -L_{\eta}H\xi +d_{*}\langle  \eta , H\xi  \rangle
+H[\xi , \eta ]_{H}-\half [H, H](\xi , \eta ))\\
&& +([\xi , \eta ]+[\xi , \eta ]_{H} ).
\ee
It thus follows that $A_{H}$ is integrable iff
for  any $\xi , \eta \in \gm (A^{*})$
\begin{equation}
\label{eq:25}
H[\xi , \eta ]=  L_{\xi}H\eta -L_{\eta}H\xi +d_{*}\langle  \eta , H\xi  \rangle
-\half [H, H](\xi , \eta ).
\end{equation}

On the other hand,
\be
(d_{*}H ) (\xi , \eta , \zeta )&=&a_{*}(\xi )\langle \eta  , H\zeta \rangle
-a_{*}(\eta )\langle \xi , H\zeta \rangle
+a_{*}(\zeta )\langle \xi , H\eta \rangle\\
 &&-\langle [\xi ,\eta ], H\zeta \rangle +\langle [\xi , \zeta ], H\eta \rangle
-\langle [\eta , \zeta ], H\xi \rangle\\
&=&\langle H[\xi , \eta  ]+L_{\eta }H\xi -L_{\xi }H\eta +d_{*}\langle \xi , H\eta \rangle , \zeta \rangle  .
\ee

 Hence,
\begin{equation}
(d_{*}H ) (\xi , \eta  )=H[\xi ,\eta ]+L_{\eta }H\xi -L_{\xi }H\eta
-d_{*}\langle \eta  , H\xi  \rangle  .
\end{equation}
This implies that Equation (\ref{eq:25}) is equivalent to
$$(d_{*}H  )(\xi , \eta  ) +\half [H, H](\xi , \eta )=0, $$
or equivalently
$$d_{*}H+\half [H, H]=0. $$ \qed
{\bf Remark} (1). Because of the symmetric   role of $A$ and $A^*$,
we have the following assertion: the graph  $A_{I}=\{ X+I X|X\in \gm ( A )\}$
of   a bundle map $I: A\lon A^{*}$
defines  a Dirac subbundle iff
$I$ is skew-symmetric and
 satisfies the following Maurer-Cartan type equation:
\begin{equation}
dI+\half [I, I]=0.
\end{equation}
(2). For the  canonical  Lie bialgebroid $(TM, T^{*}M)$ where $M$
is equipped with the zero Poisson structure, Equation  (\ref{eq:Maure}) becomes
$[H,H]=0$, which is the defining equation for a Poisson structure.  On the other
hand, if we exchange  $TM$ and $T^{*}M$, and consider the Lie bialgebroid $(T^{*}M
, TM )$, the bracket term drops out of Equation (\ref{eq:Maure}), whose solutions
correspond to a presymplectic structures. Encompassing  these two
cases into a general framework was indeed the main motivation for Courant
\cite{Courant:1990} to define and study Dirac structures.\\
(3). The Maurer-Cartan equation is a kind of
integrability equation. It is also basic in deformation theory, where it may
live on a variety of differential graded Lie algebras.  It would be interesting
to place the occurrence of this equation in our theory in a more general context.
\begin{defi}
Given a Lie bialgebroid $(A, A^{*})$, a section
$H\in \gm (\wedge^{2} A)$ is called a  {\em hamiltonian
operator }  if $A_{H}$ defines a Dirac structure.
$H$ is called a {\em strong hamiltonian operator} if
$A_{\lambda H}$ are  Dirac subbundles for all $\lambda \in \reals$.
\end{defi}

\begin{cor}
For a Lie bialgebroid $(A, A^{*})$,  $H\in \gm (\wedge^{2} A)$
is a hamiltonian operator
if Equation (\ref{eq:Maure}) holds. It  is a strong hamiltonian operator
if $d_{*}H=[H, H]=0.  $
\end{cor}

For a  hamiltonian operator $H$,  $A_{H}$ is a Dirac subbundle
which is transversal to $A$ in $A\oplus A^{*}$. Therefore, $(A, A_{H})$ is a
Lie bialgebroid according to Theorem \ref{thm:main2}.
In fact, $A_{H}$ is isomorphic to $A^{*}$,
as a vector bundle, with the anchor and bracket of its Lie algebroid
structure
given respectively  by $\hat{a}_{*}=a_{*}+a\smalcirc H$
and $[\xi , \eta \hat{]}=[\xi , \eta ]+[\xi , \eta]_{H}$, for all
$\xi , \eta \in \gm (A^{*})$. In particular, if
$H$ is a strong hamiltonian operator, one obtains
a one parameter family of Lie bialgebroids transversal
to $A$, which can be considered as a deformation
of the Lie bialgebroid $(A, A^{*})$.

Conversely, any Dirac subbundle
transversal to $A$ corresponds to  a hamiltonian
operator in an obvious way.  For example,
consider the  standard Lie bialgebra $(\mathfrak{k}, \mathfrak{b})$
arising from the Iwasawa decomposition  of $\mathfrak{k}^{\complex}$
\cite{LuW}, where $\mathfrak{k} $ is  a compact semi-simple  Lie algebra and
$\mathfrak{b}$ its corresponding dual Lie algebra.
Then any real form of $\mathfrak{k}^{\complex}$
which is transversal  to $\mathfrak{b}$ will correspond
to a hamiltonian operator (see \cite{LiuQ}   for a complete
list of such real forms for simple Lie algebras).
It is straightforward to check that such a hamiltonian operator is not strong.
  On the other hand, for the Cartan subalgebra $\mathfrak{h}$ of
$\mathfrak{k}$, every element in $\wedge^{2}h$
  gives rise to  a strong hamiltonian
  operator $H: \mathfrak{k}^{*} \lon \mathfrak{k}$,

  This  gives rise to
  a deformation of
  the standard Lie bialgebra $(\mathfrak{k}, \mathfrak{b})$ (see \cite{LS}).

These examples can be generalized to
  any gauge algebroid
  associated to
  a principal $K$-bundle.

Below, we will look at hamiltonian operators in two special
cases, each of which corresponds to some familiar
objects.

\begin{ex}
Let $P$ be a Poisson manifold with Poisson tensor $\pi$.
Let  $(TP, T^* P)$ be  the canonical
Lie bialgebroid associated to the  Poisson manifold $P$
and $E=TP \oplus T^* P$ equipped with the
induced Courant algebroid structure.
It is easy to see that   a bivector field $H$ is a hamiltonian operator
iff $H+\pi$ is a Poisson tensor. $H$ is a strong hamiltonian operator
iff $H$ is a Poisson tensor Schouten-commuting with $\pi$.
\end{ex}

\begin{ex}
Similarly, we may  switch $TP$ and $T^{*}P$, and
 consider the Lie bialgebroid $(T^{*}P, TP )$ associated to
 a Poisson manifold  $P$ with Poisson structure  $\pi$.
Let $E= T^{*}P  \oplus TP$ be equipped with its
Courant algebroid structure. In this case, a
hamiltonian operator corresponds to a two-form
$\omega \in \Omega^{2}(P)$ satisfying
$d\omega + \half [\omega , \omega ]_{\pi }=0$.
Here, $[\cdot , \cdot ]_{\pi}$ refers to the Schouten bracket
of differential forms on $P$ induced by the Poisson structure
$\pi$.
Given a hamiltonian operator $\omega$, its
graph $A_{\omega }$ defines a Dirac subbundle transversal to
$T^{*}P$, the first component of $E$ also being considered as
a Dirac subbundle. Therefore, they  form
a Lie bialgebroid. Their induced Poisson structure
on the base space can be easily checked to be given
by
$-2(\pi^{\#}+N\pi^{\#}) $, where
 $N: TP \lon TP$  is the composition
$\pi^{\#}\smalcirc  \omega^{b}$ and $\omega^{b} : TP \lon T^{*}P$
is the bundle map induced by the two-form $\omega$.
If $\omega $  is a strong hamiltonian  operator,
then $N \pi^{\#} $ defines  a Poisson structure compatible with
$\pi$. In fact, in this case, $(\pi , N)$ is
 a Poisson-Nijenhuis structure
in the sense of \cite{K-SM:1990}.
\end{ex}

We note that Vaisman \cite{va:complementary} has studied 2-forms on
Poisson manifolds satisfying the condition $[\omega , \omega
]_{\pi}=0$.    Such forms, called {\em  complementary} to the Poisson
structure, also give rise to new Lie algebroid structures on $TM$.

To end this section, we
describe a  example of Lie bialgebroids, where both  the
algebroid and its dual   arise  from
hamiltonian operators.

\begin{pro}
Let $U$ and  $V$ be Poisson tensors over  a manifold $M$
and denote by $T^{*}M_{U}$ and $T^{*}M_{V}$ their   associated
canonical cotangent Lie algebroids on  $T^{*}M$.
Assume that $U-V$ is nondegenerate. Then $(T^{*}M_{U}, T^{*}M_{V})$
is  a Lie bialgebroid, where their pairing is given
by   $(\xi , \eta )=(U-V)(\xi , \eta )$ for any
$\xi \in T^{*}M_{U}$ and
$\eta \in T^{*}M_{V}$. Their induced Poisson tensor on the
base space $M$ is given by
$-2U(U-V)^{-1}V$.
\end{pro}
\pf Let $E=TM\oplus T^* M$ be equipped with the usual
Courant bracket.
Since both $U$ and $V$ are Poisson tensors,
their graphs $A_{U}$ and $A_{V}$ are  Dirac subbundles.
They are transversal since $U-V$ is nondegenerate.
Therefore, $(A_{U}, A_{V})$ is a Lie bialgebroid,
where their pairing is given by
\begin{equation}
\llangle  U\eta +\eta , V\xi +\xi \rrangle = \half \langle \xi , (U-V)\eta \rangle .
\end{equation}
On the other hand, as Lie algebroids,
$A_{U}$ and $ A_{V}$ are clearly isomorphic to
cotangent Lie algebroids $T^{*}M_{U}$ and $ T^{*}M_{V}$
respectively.
Moreover, their anchors $a_{U}: A_{U}\lon TM $ and $a_{V}: A_{V}\lon TM $
are given respectively by
$a_{U} (U\xi +\xi )=U\xi $ and
$a_{V} (V\xi +\xi )=V\xi $.

This proves the first part of the proposition.
 To calculate their induced  Poisson structure on the base $M$,
 we need to find out
the dual map $a_{V}^{*}: T^{*}M\lon A_{V}^{*}\cong A_{U}$.
For any $\xi \in T^{*}M$,
 we assume that $a_{V}^{*} (\xi )=U\eta +\eta \in A_{U}$ via the
identification above.
For any $\zeta \in T^{*}M$,
\be
( a_{V}^{*} \xi ,  V\zeta +\zeta   )&=& \langle \xi , a_{V}
(V\zeta +\zeta )\rangle \\
&=& \langle \xi , V\zeta \rangle
\ee
On the other hand,  $( a_{V}^{*} \xi ,  V\zeta +\zeta   )=
\llangle U\eta +\eta , V\zeta +\zeta \rrangle
=\half \langle \zeta  , (U-V)\eta \rangle$.
It thus follows that $\eta =-2(U-V)^{-1}V\xi $.
Therefore, according to Proposition 3.6 in
\cite{MackenzieX:1994}, the induced
Poisson structure $a_{U}\smalcirc a_{V}^{*}:
T^{*}M\lon TM$ is given by $(a_{U}\smalcirc a_{V}^{*})\xi
=-2U(U-V)^{-1}V\xi $. \qed

Replacing $V$ by $-V$ in the proposition above, we obtain the following
``composition law'' for Poisson structures.
\begin{cor}
\label{cor-composition}
Let $U$ and  $V$ be Poisson tensors over manifold $M$
such  that $U+V$ is nondegenerate.
Then, $U(U+V)^{-1}V$ also defines a Poisson tensor on $M$.
\end{cor}

Note that, if $U$ and $V$ are nondegenerate, then
$U(U+V)^{-1}V=(U^{-1}+V^{-1})^{-1}$ is the Poisson tensor corresponding to the
sum of the symplectic forms for $U$ and $V$.  Since the sum of closed forms is
closed, it is obvious in this case that $U(U+V)^{-1}V$ is a Poisson tensor.  We do
not know such a simple proof of Corollary \ref{cor-composition} in the general
case.


\section{Null Dirac structures and Poisson reduction}

In this section, we consider another class of Dirac structures related
to Poisson reduction and dual pairs of Poisson manifolds.

\begin{pro}
\label{pro:null}
Let $(A, A^{*} )$ be a Lie bialgebroid, and  $h \subseteq A$ a subbundle of $A$.
Then $L=h\oplus h^{\perp} \subseteq  E=A\oplus A^* $
 is a Dirac structure iff  $h$   and $h^{\perp}$ are, respectively,
Lie subalgebroids of $A$ and $A^* $.
\end{pro}  
\pf  Obviously, $L=h\oplus h^{\perp}$ is a maximal isotropic subbundle of $E$.
If $L$ is  a Dirac structure, clearly  $h$ and $h^{\perp}$
are   Lie subalgebroids of $A$ and  $A^*$ respectively.
Conversely, suppose that both $h$ and $h^{\perp}$ are
Lie subalgebroids of $A$   and $A^*$ respectively. To prove
that $L$ is a Dirac structure, it suffices to show that
$[X, \xi ]$ is a section of $L$ for any
 $X \in \gm   (h) $  and $\xi \in \gm  (h^{\perp})$. Now
 $$[X,\xi]=L_{X}\xi-L_{\xi}X .$$
For any section $Y\in \gm   (h)$,
$$<L_{X}\xi, \ Y>=a(X)<\xi,Y>-<\xi,[X,Y]>=0. $$
Therefore,  $L_{X}\xi  $ is still a section of $h^{\perp}$.
Similarly, $L_{\xi}X $  is  a section of $h$.
This concludes the proof of the proposition. \qed

It is clear that a subbundle $L\subseteq E $ is of the   form $L=h\oplus h^{\perp}$
iff the minus  two-form $(\cdot,\cdot)_{-}$ on $E$,  as defined 
by Equation (\ref{eq:pairing}), vanishes  on $L$.
For this reason,  we call a Dirac structure of this
form  a  {\em  null Dirac structure}.

An immediate consequence of Proposition \ref{pro:null}
is the following:

\begin{cor}
\label{cor:reduction} 
Let $(P, \ \pi )$   be a Poisson manifold,  and $D$ a subbundle of $TP$.
Let $T^{*}P$ be equipped with the cotangent Lie
algebroid structure so that $(TP, \ T^{*}P)$ 
is  a Lie bialgebroid. Then $L=D\oplus D^{\perp} $
 is a Dirac structure  in $E= TP\oplus T^{*}P$
 iff $D$  is an integrable distribution  and
the Poisson structure on $P$ descends to 
a Poisson structure on the quotient space $P/D $ \footnote{When the
quotient  space is not  a manifold,  
this means that at each point there is a local
neighborhood $U$ such that  the Poisson structure
on $U$ descends to its quotient.} such that the natural projection
is a  Poisson map.
\end{cor} 
\pf  This follows directly  from the following 

\begin{lem}
\label{lem:reduction}
 Let  $D$  be   an integrable  distribution on a  Poisson 
manifold $P$.
$P/D$ has   an induced  Poisson structure (in the above general
sense) iff $D^{\perp}\subset T^{*}P$ is 
a subalgebroid\footnote{Such a foliation is also called
cofoliation by Vaisman \cite{vaisman:book}}.
\end{lem}
\pf For simplicity, let us assume that $P/D$ is a manifold.
The general case  will follow from the same principle.
It is clear that a function $f$  is constant along leaves of $D$ iff
  $df $  is a section of $D^{\perp}$. If $D^{\perp}$ is a
subalgebroid, it follows from the equation
\begin{equation}
\label{eq:dfg}
           d\{f,\ g\}=[df, dg] 
\end{equation}
that   $C^{\infty}(P/D)$ is a Poisson algebra.

Conversely,  a  local one-form $fdg$ is in $D^{\perp}$
iff $g$ is constant along $D$. The conclusion thus
follows from Equation (\ref{eq:dfg}) together with
the Lie algebroid axiom relating the bracket and
anchor. \qed
{\bf Remark.} \ Poisson reduction was considered by Marsden
and Ratiu in \cite{MR}.  Lemma \ref{lem:reduction} can be
considered as a special case of their theorem when $P=M$
in the Poisson triple $(P, \ M, \ E)$ (see \cite{MR}).
It would be interesting  to interpret their general
reduction theorem in terms of Dirac structures 
as in Corollary \ref{cor:reduction}.

The rest of the section is devoted to  several examples
of Corollary \ref{cor:reduction}, which will lead to some 
familiar results in Poisson geometry.

  Recall that, given a
 Poisson Lie group $G$ and a Poisson manifold $M$, an action
$$ \sigma : G\times M \lon M$$
is called a Poisson action if $\sigma$ is a Poisson map.
In this case, $M$ is called a Poisson $G$-space.

Now consider $P=G \times M$ with the product Poisson structure and
diagonal $G$-action.  Then $P/G$ is isomorphic 
to $M$, and the projection  from $P$ to $P/G=M$ 
becomes the action map $\sigma$, which is a 
 Poisson map when  $P/G$  is equipped  with the given
Poisson structure on $M$.  
  By Corollary \ref{cor:reduction}, we obtain
 a null Dirac structure $L=D\oplus D^{\perp}$ in
$TP\oplus T^*P$.

Clearly, $L$ is a Lie algebroid over $P$, which is $G$-invariant. 
It would be interesting to explore the relation between this
algebroid and  the one defined on 
$ (M\times \frakg)\oplus T^*M $, which was studied by Lu in \cite{Lu}.

For  a Poisson Lie group $G$ with  tangential Lie bialgebra
$(\frakg ,\frakg^{*})$, the Courant algebroid double $E=
 TG \oplus T^{*}G$ can be identified, as a vector
bundle, with  the trivial product
$G \times (\frakg \oplus \frakg^{*}) $ via left translation.
Under such an identification, a left invariant
 null  Dirac structure has the form $L=G\times (h\oplus h^{\perp})$,
where $h$ is  a subalgebra of $\frakg$ and $h^{\perp}$ is
a subalgebra of $\frakg^{*}$. Thus, one obtains
the following reduction 
theorem:
for a connected closed subgroup $H$ with  Lie algebra $h$, $G/H$
 has an  induced Poisson
structure  iff $h^{\perp}$  is a subalgebra of $\frakg^*$.

More generally,  let $G$ be a Poisson group,  $M$ a Poisson $G$-space.
 Suppose that $H\subseteq G$ is a closed subgroup  with Lie algebra $h$.
Assume  that $M/H$ is a nice manifold such that  the projection
$p: M\lon M/H$ is a submersion. Then the $H$-orbits define an
integrable distribution $\frakh$ on $M$. According to Corollary \ref{cor:reduction},
the Poisson structure on $M$ descends to $M/H$  
iff $\frakh^{\perp}$ is a subalgebroid 
of the cotangent algebroid $T^{*}M$ of the Poisson manifold $M$.
On the other hand,  we have 

\begin{pro}
\label{pro-perp}
If $h^{\perp}$ is a subalgebra of $\frakg^*$, then $\frakh^{\perp}$
is a subalgebroid of $T^{*}M$. Conversely, if the 
isotropic subalgebra at each point is a subalgebra of $h$,
 and in particular if the action is locally free, then
that $\frakh^{\perp}$ is a subalgebroid  implies
that $h^{\perp}\subseteq \frakg^*$ is a subalgebra.
\end{pro}
\pf It is easy to see that $\frakh^{\perp}\cong \phi^{-1}(h^{\perp})$, where
$\phi : T^{*}M\lon \frakg^{*}$ is the momentum mapping for the
lifted $G$-action on $T^{*}M$, equipped with the canonical
cotangent symplectic structure.   
 According to  Proposition 6.1 in \cite{Xu}, 
$\phi : T^{*}M\lon \frakg^{*}$ is a Lie algebroid morphism.  Before
continuing, we need the following

\begin{lem}
Let $A\lon M$ be  a Lie algebroid with anchor $a$, $\frakg$ a Lie algebra, and $\phi :A\lon
\frakg$ an algebroid morphism. Suppose that $h\subseteq \frakg$ is a subalgebra
such that $\phi^{-1}h \subseteq A$ is a subbundle. Then $\phi^{-1}h$ is  a
subalgebroid.

Conversely, given a subalgebroid $B\subseteq A$, if $\phi (B|_{m})$ is
independent of $m\in M$, then it is a subalgebra of $\frakg$.
\end{lem}
 \pf This follows directly  from the following equation
(see  \cite{Mackenzie:1992}):
$$
\phi [X,Y]=(aX)(\phi Y)-(aY)(\phi X)+{[\phi X,  \ \phi Y]}^{.}, \ \ \forall
X, Y\in \gm (A),
$$
where $\phi X$,   $\phi Y$ and $\phi [X, Y]$ are considered as
$\frakg$-valued functions on $P$, and $[\cdot , \cdot]^{\cdot}$
refers to  the pointwise bracket. \qed

Now, the first part of Proposition \ref{pro-perp}
 is obvious according to the lemma above. For the
second part, we only need to note that 
 $\phi (\frakh^{\perp} ) =  h^{\perp}\cap Im \phi $, and the assumption that
the isotropic subalgebra at each  point is a subalgebra of $h$ is equivalent
to that $ h^{\perp}\subseteq  Im \phi $.
This concludes the proof of the Proposition.  \qed

   From the   above discussion, we   have the following conclusion:
if $h^{\perp}$ is a subalgebra of $\frakg^*$, then  the Poisson structure
on $M$ descends to $M/H$. This is a well-known reduction theorem of
Semenov-Tian-Shansky \cite{STS} (see also \cite{we:coisotropic}).

Conversely,  if  the 
isotropic subalgebra  at each point is a subalgebra of $h$,
 and in particular  if  the action is locally free,  the converse
also holds.

Another interesting example arises when   $P$ is  a symplectic manifold
 with an  invertible Poisson tensor $\pi$.
     In this case, $\pi^{\#}: \ T^{*}P\lon TP$ is a Lie algebroid isomorphism.
Given  a null Dirac structure  $L=D\oplus D^{\perp}$,
$\bar{D}= \pi^{\#}(D^{\perp})$ is   a subalgebroid of
$TP$. It is simple to see that   $(\bar{D})^{\perp}=(\pi^{\#})^{-1}(D)$, and
is therefore  a subalgebroid of $T^{*}P$.
 Thus, $\bar{L} \stackrel{def}{=}\bar{D}\oplus (\bar{D})^{\perp}$
defines another null Dirac structure.
 It is easy to see that $D$ and $\bar{D}$ are symplectically
orthogonal to each another. Thus $P/\bar{D}$ is  a Poisson
manifold (assume that it is a nice manifold) so that
$P/D$ and $P/\bar{D}$  constitute a   full dual pair,
which  is a well known result of Weinstein \cite{we:1983}.
Conversely, it is clear that  a full dual pair  corresponds to  
  a null Dirac structure.


\begin{thebibliography}{99}
\bibitem{ba:quasi}
Bangoura, M., Quasi-groupes de Lie-Poisson,
{\em C.R. Acad. Sci. Paris} {\bf 319} (1994), 974-978.
\bibitem{ba-ko:double}
Bangoura, M. and Kosmann-Schwarzbach, Y.
     The double of a Jacobian quasi-bialgebra, {\em Lett. Math. Phys.}
     {\bf 28} (1993), 13-29.
\bibitem{br-da-ha:topological}
Brown, R., Danesh-Naruie, G., and Hardy, G.P.L., Topological
groupoids: II. Covering morphisms and $G$-spaces,
{\em Math. Nachr.} {\bf 74} (1976), 143-156.
\bibitem{Courant:1990}
Courant, T.J., Dirac manifolds, {\em Trans. A.M.S.} {\bf 319} (1990),
631-661.
\bibitem{Dorfman}
Dorfman, I.Ya., {\em Dirac structures and integrability of nonlinear evolution
  equations}, Wiley, Chichester, 1993.
\bibitem{dr:quantum}
Drinfel'd, V.G., Quantum groups, {\em Proc. ICM, Berkeley,} 1986, 
789-820.
\bibitem{dr:quasi}
Drinfel'd, V. G., Quasi-Hopf algebras, {\em Leningrad Math. J.} {\bf 2}
(1991), 829-860.
\bibitem{dr:poisson}
Drinfel'd, V.G., On Poisson homogeneous spaces of Poisson-Lie groups,
{\em Theor. Math. Phys.} {\bf 95} (1993), 524-525.
\bibitem{gi:resolution}
Ginzburg, V.A., Resolution of diagonals and moduli spaces, {The Moduli
Space of Curves}, R. Dijkgraff et al, eds., Birkh\"auser Boston, (1995),
231-266.
\bibitem{hu:poisson} 
Huebschmann, J., Poisson cohomology and quantization,
{\em J. Reine Angew. Math.} {\bf 408} (1990), 57-113.
\bibitem{ko:quasi}
Kosmann-Schwarzbach, Y., Jacobian quasi-bialgebras and quasi-Poisson
Lie groups, {\em Contemp. Math.} {\bf 132} (1992), 459-489.
\bibitem{K-S:1994}
Kosmann-Schwarzbach, Y., Exact Gerstenhaber
algebras and Lie bialgebroids, {\em Acta Appl. Math.} {\bf 41}
(1995), 153-165.
\bibitem{K-SM:1990}
Kosmann-Schwarzbach, Y., and Magri, F., Poisson-Nijenhuis structures,
{\em Ann. Inst. H.~Poincar{\'e} Phys. Th{\'e}or.} {\bf 53}, (1990),  35--81.
\bibitem{la-ma:strongly} 
Lada, T., and Markl, M., Strongly homotopy Lie algebras,
{\em Comm. in Alg.} {\bf 23} (1995), 2147-2161.
\bibitem{LS} Levendorskii, S., and Soibelman, Y., Algebras of functions on
compact quantum groups, Schubert cells and quantum tori,
{\em Comm. Math. Phys.} {\bf  139 } (1991), 141-170.
\bibitem{LiuQ} Liu, Z.-J. and Qian, M., Generalized Yang-Baxter
equations, Koszul operators and Poisson Lie groups, {\em J. Diff.
Geom.} {\bf 35} (1992), 399-414.
\bibitem{LWX}
Liu, Z.-J., Weinstein, A.,  and Xu, P.,
Dirac structures and Poisson homogeneous spaces,
preprint, 1996.
\bibitem{Lu} Lu, J.-H., Lie algebroids associated to Poisson actions,
{\em Duke Math. J.} (to appear).
\bibitem{LuW} Lu, J.-H. and Weinstein, A.,   Poisson Lie groups,
dressing transformations, and the Bruhat decomposition, {\em J. Diff.
Geom.} {\bf 31} (1990), 501-526.
\bibitem{Mackenzie:book}
Mackenzie, K., {\em Lie Groupoids and Lie Algebroids in Differential
Geometry}, LMS Lecture Notes Series, {\bf 124}, Cambridge Univ. Press, 1987.
\bibitem{Mackenzie:1992}
Mackenzie, K.,  Double Lie algebroids and second-order geometry I,
{\em Adv. in Math.} {\bf 94} (1992),  180-239. 
\bibitem{MackenzieX:1994}
Mackenzie, K.C.H. and Xu, P.,
Lie bialgebroids and Poisson groupoids,
{\em Duke Math. J.} {\bf 18} (1994), 415-452.
\bibitem{MackenzieX:1996}
Mackenzie, K.C.H. and Xu, P.,
Integration of Lie bialgebroids,  preprint, 1996.
\bibitem{MR}
Marsden, J.E., and Ratiu, T., Reduction of Poisson manifolds,
{\em Lett. Math. Phys.} {\bf 11} (1986), 161-169.
\bibitem{me-re:algebres}
Medina, A., and Revoy, P., Alg\`{e}bres de Lie et produit scalaire
invariant, {\em Ann. Sc. Ec. Norm. Sup., $4^{e}$ s\'{e}rie} {\bf 18}
(1985), 553-561.
\bibitem{me-re:lie}
Medina, A., and Revoy, P., Groupes de Lie-Poisson et double extension,
{\em Seminaire Gaston-Darboux de G\'{e}om\'{e}trie et Topologie
Diff\'{e}rentielle 1990-1991, Universit\'{e} Montpellier II} (1992), 87-105.
\bibitem{mi-we:moments}
Mikami, K., and Weinstein, A., Moments and reduction for symplectic
groupoid actions, {\em Publ. RIMS Kyoto Univ.} {\bf 24}
(1988),121-140.
\bibitem{mi:rational}
Millson, J., Rational homotopy theory and deformation problems from
algebraic geometry, {\em Proc. I.C.M., Kyoto, Japan, 1990},
Springer-Verlag, Tokyo, 1991, 549-558.
\bibitem{STS}
Semenov-Tian-Shansky, M.A., Dressing transformations and Poisson
group actions, 
{\em Publ. RIMS, Kyoto University}
{\bf 21} (1985), 1237-1260.
\bibitem{vaisman:book}
Vaisman, I., {\em Lectures on the geometry of Poisson manifolds},
PM {\bf 118}, Basel; Boston; Berlin; Birkh\"auser 1994.
\bibitem{va:complementary}
Vaisman, I., Complementary 2-forms of Poisson structures, {\em Compositio
Math.}, to appear.
\bibitem{we:1983}
Weinstein, A., The local structure of Poisson manifolds,
{\em J. Diff. Geom.} {\bf 18} (1983), 523-557.
\bibitem{we:coisotropic}
Weinstein, A., Coisotropic calculus and Poisson groupoids, {\em J.
Math. Soc. Japan} {\bf 40} (1988), 705-727.
\bibitem{Xu}
Xu, P., On Poisson groupoids, {\em International J. of Math.}
{\bf 6} (1995), 101-124.

\end{thebibliography}
     \end{document}